\documentclass[graybox]{svmult}
\usepackage{geometry}                % See geometry.pdf to learn the layout options. There are lots.
\usepackage{graphicx}
\usepackage{amssymb}
\usepackage{amsmath} %loaded before
\usepackage{color}

\usepackage{bm}
\usepackage{epstopdf}
\usepackage{accents} %loaded before
\usepackage{stmaryrd}
\usepackage{adjustbox}

\definecolor{redcol}{rgb}{1.,0.,0.0} 
\definecolor{commcol}{rgb}{0.,0.,1.}
\definecolor{lnkcol}{rgb}{0.,0.,0.0} %black
\usepackage[pdftex,pdfpagelabels=true,bookmarksdepth=3]{hyperref}
\hypersetup{  colorlinks=true,
              citecolor=lnkcol,%
              linkcolor=lnkcol,%
              urlcolor=lnkcol,
            pdfborder=false,
            bookmarks=true,
            bookmarksnumbered=true,
            pdftitle = {Stability of Wall Boundary Condition Procedures for Discontinuous Galerkin Spectral Element Approximations of the Compressible Euler Equations},
pdfsubject = {},
pdfauthor = {Florian J. Hindenlang, Gregor J. Gassner, David A. Kopriva},
pdfkeywords = {}
 }
\usepackage[all]{hypcap}  % hyperlink soll oberhalb des Bilds anhalten
\usepackage[dvips]{thumbpdf}

\usepackage{cleveref} % === > !!! MUST BE USED AFTER hyperref!!! <====

  \crefname{chapter}{Chap.}{Chaps.}       
  \Crefname{chapter}{Chapter}{Chapters}
  \crefname{section}{Sec.}{Secs.}
  \Crefname{section}{Section}{Sections}
  \crefname{table}{Table}{Tables}
  \Crefname{table}{Table}{Tables}
  \crefname{figure}{Fig.}{Figs.}
  \Crefname{figure}{Figure}{Figures}
  \crefname{equation}{}{}
  \Crefname{equation}{Equation}{Equations}
  \crefname{algorithm}{Alg.}{Algs.}
  \Crefname{algorithm}{Algorithm}{Algorithms}
  \crefname{thm}{Thm.}{Thms.}
  \Crefname{thm}{Theorem}{Theorems}
  \crefname{lem}{Lemma}{Lemmas}
  \Crefname{lem}{Lemma}{Lemmas}
  \crefname{cor}{Corr.}{Corrs.}
  \Crefname{cor}{Corrollary}{Corrollaries}
  \crefname{prop}{Prop.}{Props.}
  \Crefname{prop}{Proposition}{Propositions}
  \crefname{rem}{Rem.}{Rems.}
  \Crefname{rem}{Remark}{Remarks}
  \crefname{appendix}{Appendix}{Appendices}
  \Crefname{appendix}{Appendix}{Appendices}
  
\usepackage{stackengine}
\stackMath

\usepackage{accsupp} % for ensuring the right Unicode codepoint upon pasting

\makeatletter
\def\munderbar#1{\underline{\sbox\tw@{$#1$}\dp\tw@\z@\box\tw@}}
\makeatother

\newcommand\iprod[1]{\left\langle #1\right\rangle} 				% inner product
\newcommand\inorm[1]{\left |\left| #1\right|\right|}		% norm
\DeclareMathAccent{\spacevec}{\mathord}{letters}{126}           % same but without using the accents package, looks better, can be used in figure captions...
\newcommand\acclrvec[1]{\accentset{\,\leftrightarrow}{#1}}	% define leftrightarrow as an accent

\newcommand\contravec[1]{\tilde{ #1}}					% contravariant vector
\newcommand\statevec[1]{\mathbf #1}					% state vector, e.g. [rho, rhov, E, B]^T
\newcommand\contrastatevec[1]{\tilde{\mathbf #1}} 			% contravariant state vector
\newcommand\bigstatevec[1]{\acclrvec{{\mathbf #1}}}		% block vector
\newcommand\bigcontravec[1]{\acclrvec{\tilde{\mathbf #1}}} 	% contravariant block vector
\newcommand\oneHalf{\frac{1}{2}}

\newcommand\dS{\,\operatorname{dS} }

\newcommand\ext{\mathrm{ext}}

\newcommand\ent{{\,\varsigma}} %superscript for entropy flux

\newcommand\mmatrix[1]{\underbar{#1}}				% Matrix as taught in linear algebra, i.e. math matrix.
\newcommand\mmatrixAn{{\mmatrix{A}}_{\,n}}
\newcommand{\jump}[1]{\left\llbracket#1\right\rrbracket}   % jump at element boundary
\newcommand\Ma{{\operatorname{Ma}}}

\usepackage{type1cm}        % activate if the above 3 fonts are
\usepackage{makeidx}         % allows index generation
\usepackage{graphicx}        % standard LaTeX graphics tool
\usepackage{multicol}        % used for the two-column index
\usepackage[bottom]{footmisc}% places footnotes at page bottom

\usepackage{newtxtext}       % 
\usepackage{newtxmath}       % selects Times Roman as basic font

\makeindex             % used for the subject index
\usepackage{amssymb}
\usepackage{epstopdf}
\title*{Stability of Wall Boundary Condition Procedures for Discontinuous Galerkin Spectral Element Approximations of the Compressible Euler Equations}
\titlerunning{Stability of Wall Boundary Condition Procedures}
\author{Florian J. Hindenlang, Gregor J. Gassner and David A. Kopriva}
\institute{Florian J. Hindenlang \at  Max Planck Institute for Plasma Physics, Boltzmannstra{\ss}e 2, D-85748 Garching, Germany, \email{florian.hindenlang@ipp.mpg.de}
\and Gregor J. Gassner\at Department for Mathematics and Computer Science/Center for Data and Simulation Science, University of Cologne, Cologne, Germany \email{ggassner@math.uni-koeln.de}
\and David A. Kopriva \at Florida State University and San Diego State University \email{kopriva@math.fsu.edu}}
\begin{document}
\maketitle

\abstract{ We perform a linear and entropy stability analysis for wall boundary condition procedures for discontinuous Galerkin spectral element  approximations of the compressible Euler equations. Two types of boundary procedures are examined. The first defines a special wall boundary flux that incorporates the boundary condition. The other is the commonly used reflection condition where an external state is specified that has an equal and opposite normal velocity. The internal and external states are then combined through an approximate Riemann solver to weakly impose the boundary condition. We show that with the exact upwind and Lax-Friedrichs solvers the approximations are energy dissipative, with the amount of dissipation proportional to the square of the normal Mach number. Standard approximate Riemann solvers, namely Lax-Friedrichs, HLL, HLLC are entropy stable. The Roe flux is entropy stable under certain conditions. An entropy conserving flux with an entropy stable dissipation term (EC-ES) is also presented. The analysis gives insight into why these boundary conditions are robust in that they introduce large amounts of energy or entropy dissipation when the boundary condition is not accurately satisfied, e.g. due to an impulsive start or under resolution.}

\section{Introduction}

The ingredients for a reliable numerical method for the approximation of partial differential equations, e.g. one that will not blow up, include stable inter-element and physical boundary condition implementations. The recognition that the discontinuous Galerkin spectral element method (DGSEM) with Gauss-Lobatto quadratures satisfies a summation-by-parts (SBP) operators \cite{Gassner:2016ye,gassner2010} has allowed for the analysis of these schemes and to connect them with penalty collocation and SBP finite difference schemes. For instance, in \cite{Gassner2018}, we showed that a split form approximation of the compressible Navier-Stokes equations was both linearly and entropy stable provided that the boundary conditions were properly imposed. 

The importance of stable boundary condition procedures for hyperbolic equations has long been studied, especially in relation to finite difference methods, e.g. \cite{skew_sbp2,nordstroem:2006,Nordstrom:2016jk}. Only recently have they been studied for discontinuous Galerkin approximations. In \cite{PARSANI201588}, the authors showed that the reflection approach is stable when using an entropy conserving flux and an additional entropy stable dissipation term (EC-ES). In \cite{CHEN2017427}, the authors show that the reflection condition is stable if the numerical flux is either the Godunov or HLL flux.

In this paper, we analyze both the linear and entropy stability of two types of commonly used wall boundary condition procedures used with the DGSEM applied to the compressible Euler equations. In both cases, wall boundary conditions are implemented through a numerical flux. The boundary condition might be implemented through a special wall numerical flux that includes the boundary condition, or a fictitious external state applied to a Riemann solver approximation. We show how to construct special wall numerical fluxes that are stable, and study the behavior of the approximations. In particular, we show that the use of Riemann solvers at the boundaries introduce numerical dissipation in an amount that depends on the size of the normal Mach number at the wall.

\section{The compressible Euler Equations and the Wall Boundary Condition}
We write the Euler equations as
\begin{equation}
\label{eq:euler}
{\statevec u_t} + \sum\limits_{i = 1}^3 {\frac{{\partial {\statevec f_{i}}}}{{\partial {x_i}}}} = 0.
\end{equation}
The state vector contains the conservative variables
\begin{equation}\statevec u = \left[ {
  \varrho \;\;
  {\varrho \spacevec v} \;\;
  {E}
} \right]^{T} = \left[ {
  \varrho  \;\;
  {\varrho v_1} \;\;
  {\varrho v_2} \;\; 
  {\varrho v_3} \;\;
  {E} 
} \right]^{T}.\end{equation}
In standard form, the components of the advective fluxes are
\begin{equation}
\statevec f_{1}  = \left[ {\begin{array}{*{20}c}
   {\varrho v_1}  \\
   {\varrho v_1^2 + p}  \\
   {\varrho v_1\,v_2}  \\
   {\varrho v_1\,v_3}  \\
   {(E + p)v_1}  \\

 \end{array} } \right]\quad \statevec f_{2}  = \left[ {\begin{array}{*{20}c}
   {\varrho v_2}  \\
   {\varrho v_2\,v_1}  \\
   {\varrho v_2^2  + p}  \\
   {\varrho v_2\,v_3}  \\
   {(E + p)v_2}  \\

 \end{array} } \right]\quad \statevec f_{3}  = \left[ {\begin{array}{*{20}c}
   {\varrho v_3}  \\
   {\varrho v_3\,v_1}  \\
   {\varrho v_3\,v_2}  \\
   {\varrho v_3^2  + p}  \\
   {(E + p) v_3}  \\

 \end{array} } \right],
\end{equation}
Here, $\varrho,\,\spacevec{v}=(v_1,v_2,v_3)^T,\,p,\,E$ are the mass density, fluid velocities, pressure and total energy. We close the system with the ideal gas assumption, which relates the total energy  and pressure 
\begin{equation}
p = (\gamma-1)\left(E - \frac{1}{2}\varrho\left\|\spacevec{v}\right\|^2\right),
\label{eqofstate}
\end{equation}
where $\gamma$ denotes the adiabatic coefficient. 
For a compact notation that simplifies the analysis, we define \emph{block vectors} (with the double arrow)
\begin{equation}
\bigstatevec{f} =
 \left[ {
  {{\statevec f_1}} \;\;
  {{\statevec f_2}} \;\;
  {{\statevec f_3}} 
} 
\right]^{T}\,,
\end{equation}
so that the system of equations can be written in the compact form
\begin{equation}
\statevec u_{t} + \spacevec \nabla_x\cdot\bigstatevec{f}=0.
\label{eq:ConservationLaw}
\end{equation}

The linear Euler equations are derived by linearizing about a constant mean state $(\bar{\varrho},\bar{v}_1,\bar{v}_2,\bar{v}_3,\bar{p})$. 
We follow \cite{Nordstrom:2005qy} for the symmetrization of the linearized equations, with the constants 
\begin{equation}
a = \sqrt {\frac{{\gamma  - 1}}{\gamma }} \bar c,\quad b = \frac{\bar c}{{\sqrt \gamma  }}\,,\bar c=\sqrt{\frac{\gamma\bar{p}}{\bar{\varrho}}},
\end{equation}
where $\bar c$ is the sound speed of the constant mean state. 
The state variables become 
\begin{equation}
\statevec u = \left[\varrho'\;\; v_{1}\;\; v_{2}\;\; v_{3}\;\; p'\right]^{T},
\end{equation}
where $\spacevec{v}$ is the velocity perturbation from the mean state, and we introduce
\begin{equation}
 \varrho' = b\frac{\tilde{\varrho}}{\bar{\varrho}}\,,\quad p'=\frac{1}{\bar{\varrho}a}\tilde{p} -\frac{1}{\sqrt{\gamma-1}}\varrho'   ,
 \end{equation} 
which depend on the density and pressure perturbations $\tilde{\varrho},\tilde{p}$.
The flux vectors are
\begin{equation}
\statevec f_{i}=\mmatrix A_{i}\statevec u,\quad \bigstatevec f = \spacevec {\mmatrix A}\statevec u = \left(\mmatrix A_{1}\hat x + \mmatrix A_{2}\hat y+\mmatrix A_{3}\hat z\right)\statevec u,
\end{equation}
where \cite{Nordstrom:2005qy}
\begin{equation}{\mmatrix A_1} = \left[ {\begin{array}{*{20}{c}}
 \bar v_1&b&0&0&0 \\ 
  b&\bar v_1&0&0&a \\ 
  0&0&\bar v_1&0&0 \\ 
  0&0&0&\bar v_1&0 \\ 
  0&a&0&0&\bar v_1 
\end{array}} \right],\quad {\mmatrix A_2} = \left[ {\begin{array}{*{20}{c}}
  \bar v_2&0&b&0&0 \\ 
  0&\bar v_2&0&0&0 \\ 
  b&0&\bar v_2&0&a \\ 
  0&0&0&\bar v_2&0 \\ 
  0&0&a&0&\bar v_2 
\end{array}} \right],\quad {\mmatrix A_3} = \left[ {\begin{array}{*{20}{c}}
  \bar v_3&0&0&b&0 \\ 
  0&\bar v_3&0&0&0 \\ 
  0&0&\bar v_3&0&0 \\ 
  b&0&0&\bar v_3&a \\ 
  0&0&0&a&\bar v_3 
\end{array}} \right]\end{equation}
are constant symmetric matrices. 

The linear equations have the property that the $L^{2}$ norm of the solution over a domain $\Omega$ is bounded by terms of the boundary data on $\partial\Omega$, only. Let
\begin{equation}
\iprod{\statevec v,\statevec w} = \int\limits_{\Omega}{\statevec v^{T}\statevec w\,dxdydz },\quad \iprod {\bigstatevec f,\bigstatevec g } = \int\limits_{\Omega} {\sum\limits_{i = 1}^3 {\statevec f_i^T{\statevec g_i}} \,dxdydz }.
\end{equation}
represent the $L^{2}$ inner product of two state vectors $\statevec v$ and $\statevec w$ and two block vectors $\bigstatevec f$ and $\bigstatevec g$, respectively. Since the coefficient matrices are constant the product rule and symmetry of $\spacevec {\mmatrix A}$ implies
\begin{equation}
\iprod{\spacevec\nabla_x\cdot\bigstatevec f,\statevec u} = \iprod{\spacevec\nabla_x\cdot\left(\spacevec{\mmatrix A}\statevec u\right),\statevec u}=\iprod{\nabla_x\statevec u,\bigstatevec f}.
\end{equation}
Then it follows from Gauss' law (integration by parts) that 
\begin{equation}
\iprod{\spacevec\nabla_x\cdot\bigstatevec f,\statevec u}=\oneHalf\int_{\partial\Omega} \statevec  u^{T}\bigstatevec f\cdot\spacevec ndS,
\end{equation}
where $\spacevec n$ is the outward normal to the surface of $\Omega$.
The norm of the solution therefore satisfies
\begin{equation}
\frac{d}{dt}\inorm{\statevec u}^{2}=-\int_{\partial\Omega} \statevec u^{T}\bigstatevec f\cdot\spacevec ndS.
\label{eq:ContinuousEnergy}
\end{equation}

Replacing the boundary terms by boundary conditions leads to a bound on the solution in terms of the boundary data. The argument of the boundary integral on the right of \eqref{eq:ContinuousEnergy} is
\begin{equation}
\statevec u^{T}\bigstatevec f\cdot\spacevec n=\statevec u^{T}\left(\spacevec{\mmatrix A}\cdot\spacevec n\right)\statevec u = 2\left(\varrho b + ap\right)v_{n}+(\spacevec{\bar{v}}\cdot\spacevec{n})(\varrho^2 + |\spacevec{v}|^2+p^2),
\end{equation}
where $v_{n}$ is the wall normal velocity, $v_{n}=\spacevec v \cdot\spacevec n$. 
Note that here, the mean flow must be chosen such that the normal flow vanishes at the wall boundary $\spacevec{\bar{v}}\cdot\spacevec{n}=0$, so that the boundary condition makes physical sense.

Therefore, with the no penetration wall condition $v_{n}=0$ applied,
\begin{equation}
\frac{d}{dt}\inorm{\statevec u}^{2} = 0,
\end{equation}
and the (energy) norm of the solution is bounded for all time by its initial value.

The nonlinear equations, on the other hand, satisfy a bound on the entropy that depends only on the boundary data. For what follows, we assume that the solution is smooth so that we don't have to consider entropy generated at shock waves. We introduce the entropy density (scaled with $(\gamma-1)$ for convenience) as
\begin{equation}
 s(\statevec u) = -\frac{\varrho \varsigma}{(\gamma-1)}\,,
\end{equation} 
where $\varsigma=\ln(p)-\gamma\ln(\varrho)$ is the physical entropy. (The minus sign is
conventional in the theory of hyperbolic conservation laws to ensure a decreasing entropy function.)
The entropy flux for the Euler equations is
\begin{equation}
 \spacevec f^\ent(\statevec u) = \spacevec v \, s = -\frac{\varrho \varsigma \spacevec v}{(\gamma-1)}\,.
\end{equation} 
Finally the entropy variables are
\begin{equation}
 \statevec w = \frac{\partial s(\statevec u)}{\partial\statevec u}= \left[ {\begin{array}{*{20}{c}}
  \frac{\gamma-\varsigma}{\gamma-1} - \beta ||\spacevec v ||^2,  \\ 
  {2 \beta \spacevec v} \\ 
  {-2\beta} 
\end{array}} \right]\,,\quad \beta = \frac{\varrho}{2p}\,.
\end{equation} 
The entropy pair contracts the solution and fluxes, meaning that it satisfies the relations
\begin{equation}
\statevec{w}^T\,\statevec{u}_t = \left(\frac{\partial s}{\partial\statevec{u}}\right)^T\statevec{u}_t  = s_t(\statevec{u}),\quad \statevec{w}^T\, \spacevec{\nabla}_x\cdot\bigstatevec{f} = \spacevec\nabla_x\cdot\spacevec{f}^\ent.
\label{eq:wsContraction}
\end{equation}

When we multiply \eqref{eq:ConservationLaw} with the entropy variables and integrate over the domain,
\begin{equation}
\iprod{\statevec{w}(\statevec u), {\statevec u_t}} + \iprod{\statevec{w}(\statevec u),\spacevec\nabla_x\cdot\bigstatevec{f}} = 0\,. \label{eq:weakformContracted}
\end{equation}
Next we use the properties of the entropy pair to contract \cref{eq:weakformContracted} and use integration by parts to get
\begin{equation}
\iprod{s_t(\statevec u ),1} =- \iprod{\spacevec\nabla_x\cdot\spacevec{f}^\ent,1} = -\int\limits_{\partial\Omega}\left(\spacevec{f}^\ent\cdot\spacevec{n}\right) \dS\,
\label{eq:weakformContract2}
\end{equation}
showing that, in the continuous case, the total entropy in the domain can only change via the boundary conditions. 

In the case of a zero-mass flux boundary condition, with $v_n=\spacevec v \cdot \spacevec n=0$, the entropy is not changed by the slip-wall boundary condition, since 
\begin{equation}
 -\spacevec f^\ent\cdot \spacevec n=\frac{\varrho \varsigma}{(\gamma-1)} v_n=0.
\end{equation} 

\section{Stability bounds for the DGSEM}
The DGSEM is described in detail in \cite{Gassner2018} and elsewhere \cite{Black1999,Kopriva:2009nx}. We will only quickly summarize the approximation here. The domain, $\Omega$ is subdivided into non-overlapping, conforming, hexahedral elements. Each element is mapped to the reference element $E=[-1,1]^{3}$. Associated with the transformation from the reference element is a set of contravariant coordinate vectors, $\spacevec a^{i}$, and transformation Jacobian, $\mathcal J$. The equations \eqref{eq:ConservationLaw} transform to another conservation law on the reference element as
 \begin{equation}
 \mathcal J {{\statevec u}_t} + {\spacevec\nabla _\xi } \cdot {\bigcontravec f} = 0,
\label{eq:xFormedEulerEquations}
\end{equation}
where $\bigcontravec f$ is the contravariant flux vector with components $\contravec{\statevec f}^{i} = \mathcal J \spacevec a^{i}\cdot\bigstatevec f$. 

The approximation of \eqref{eq:xFormedEulerEquations} proceeds as follows: A weak form is created by taking the inner product of the equation with a test function. The Gauss law is applied to the divergence term to separate the boundary from the interior contributions. The resulting weak form is then approximated:  The solution vector is approximated by a polynomial of degree $N$ interpolated at the Legendre Gauss-Lobatto points. In the following, we will represent the true continuous solutions by lower case letter. Upper case letters will denote their polynomial approximations, except for the density, where the approximation is denoted by $\rho$. The volume fluxes are replaced by two-point numerical fluxes. In the linear case, the two point fluxes are immediately relatable to a split form of the equations. Integrals are replaced by Legendre-Gauss-Lobatto quadratures. Finally, the boundary fluxes are replaced by a numerical flux. See \cite{Gassner2018} and \cite{Gassner:2013uq} for details.

The result is an approximation that is energy stable for the linearized equations if at every quadrature point along a physical boundary the numerical flux $\contrastatevec F^{*}$ satisfies the bound \cite{Gassner2018}
\begin{equation}
{\statevec U^{T}{\left\{{\contrastatevec F}^{*}-\frac{1}{2}\bigcontravec F\cdot\hat n\right\}}}\ge 0,
\label{eq:LinearStabilityConditionCF}
\end{equation}
where $\bigcontravec F$ is the polynomial interpolation of the contravariant flux from the interior, $\hat n$ is the reference space outward normal direction, and $\statevec U $ is the approximation of the state vector.
Since the contravariant fluxes are proportional to the normal fluxes \cite{Kopriva:2009nx}, we can change the condition \eqref{eq:LinearStabilityConditionCF} to
\begin{equation}
B_{L}\equiv{\statevec U^{T}{\left\{{\statevec F}^{*}-\frac{1}{2}\bigstatevec F\cdot\spacevec n\right\}}}\ge 0,
\label{eq:LinearStabilityCondition}
\end{equation}

For entropy stability of the nonlinear equations, the boundary stability condition shown in  \cite{Gassner2018} is proportional to
\begin{equation}
B_{NL} \equiv  \statevec{W}^T\left(\statevec{F}^{*}-\left(\bigstatevec{F}\cdot\spacevec{n}\right)\right) + \left(\spacevec{F}^\ent\cdot\spacevec{n}\right)\ge 0,
\label{eq:EntropyStabilityCondition}
\end{equation}
where $\spacevec{F}^\ent$ is the polynomial interpolation of the entropy flux, $\spacevec f^{\ent}$, and $\statevec W$ is the interpolation of the entropy variables.

\subsection{Linear Stability of Wall Boundary Condition Approximations}

To find linearly stable implementations of the wall condition $v_{n}=0$, one needs only find a numerical flux that satisfies it and the condition \eqref{eq:LinearStabilityCondition}. For the linear equations, the approximation of the state vector is $\statevec U= [\rho'\;\; \spacevec V \;\; P']^T$ and the normal contravariant flux is proportional to
\begin{equation}
\bigstatevec F\cdot\spacevec n = \spacevec{ {\mmatrix A}}\cdot\spacevec n \,\statevec U = \left[ bV_{n}\;\; n_1 Q\;\; n_2 Q\;\;n_3 Q\;\; aV_{n} \right]^{T},
\label{eq:InteriorNormalFlux}
\end{equation}
where $V_{n}$ is the approximation of the normal velocity at the wall computed from the interior, $Q = b\rho' + aP'$, and $(n_1,n_1,n_3)$ are the three components of the physical space normal vector, $\spacevec n$. The numerical flux can be expressed as
\begin{equation}
\statevec{F}^{*} = \spacevec{ {\mmatrix A}}\cdot\spacevec n \,\statevec U^{*}= \left[ bV^{*}_{n}\;\; n_1 Q^{*}\;\; n_2 Q^{*}\;\;n_3 Q^{*}\;\; aV^{*}_{n} \right]^{T}.
\label{eq:NormalStarFlux}
\end{equation} It then remains only to find $Q^{*}$ so that \eqref{eq:LinearStabilityCondition} is satisfied when the normal wall condition $V^{*}_{n}=0$ is applied. When we substitute the fluxes \eqref{eq:InteriorNormalFlux} and \eqref{eq:NormalStarFlux} into \eqref{eq:LinearStabilityCondition}, 
\begin{equation}
B_{L}=\oneHalf\left\{ Q\left( 2V^{*}_{n}-V_{n}\right) + V_{n}\left( 2Q^{*}-Q\right)\right\} = \oneHalf\left\{ 2QV^{*}_{n} + 2V_{n}\left( Q^{*}-Q\right)\right\}
\label{eq:DirectNumericalFlux}
\end{equation}
Substituting the wall boundary condition $V^{*}_{n}=0$ yields the condition on $Q^{*}$ for stability
\begin{equation}
V_{n}\left( Q^{*}-Q\right)\ge 0.
\end{equation}
Neutral stability is thus ensured if $\rho^{*}$ and $P^{*}$ are computed from the interior, i.e. $\rho^{*}=\rho'$, $P^{*}=P'$ so that $Q^{*}=Q$.

In practice, the boundary condition is also implemented through the use of a Riemann solver and external state designed to imply the physical boundary condition to construct the numerical boundary flux. The exact upwind ($\varepsilon=1$) normal Riemann flux and the central flux ($\varepsilon=0$) for the linear system of equations is
\begin{equation}
\statevec F^{*}\left( \statevec U, \statevec U^{\ext}\right) = \frac{1}{2}\left\{\bigstatevec F\left(\statevec U\right)\cdot\spacevec n + \bigstatevec F\left(\statevec U^{\ext}\right)\cdot\spacevec n\right\} - \frac{\varepsilon}{2}\left|{\mmatrixAn}\right|\left(\statevec U^{\ext} - \statevec U\right),
\end{equation}
where ${\mmatrixAn}\equiv {\spacevec{\mmatrix A}}\cdot\spacevec n$ is the normal coefficient matrix. The external state is set by by using the interior values of the density and pressure and the negative of the value of the normal velocity, 
\begin{equation}
\statevec U^{\ext}=\left[ {
  \rho'  \;\;
  {} \;\;
  {\left(\spacevec V - 2V_{n}\spacevec n\right)} \\ 
  {} \;\;
  P' 
} \right]^{T}.
\label{eq:LinearExtState}
\end{equation}

For $\varepsilon=0$, using the central (averaged) numerical flux, the interior flux contribution cancels and condition \eqref{eq:LinearStabilityCondition} reduces to
\begin{equation}
B_{L,0}=\oneHalf\statevec U^{T} {\mmatrixAn}\statevec U^{\ext} = \left[ {\begin{array}{*{20}{c}}
  \rho' &{}&\spacevec V&{}&P' 
\end{array}} \right]\left[ {\begin{array}{*{20}{c}}
  0&{n_1 b}&{n_2 b}&{n_3 b}&0 \\ 
  {n_1 b}&0&0&0&{n_1 a} \\ 
  {n_2 b}&0&0&0&{n_2 a} \\ 
  {n_3 b}&0&0&0&{n_3 a} \\ 
  0&{n_1 a}&{n_2 a}&{n_3 a}&0 
\end{array}} \right]\left[ {\begin{array}{*{20}{c}}
  \rho'  \\ 
  {} \\ 
   {\spacevec V - 2V_n\spacevec n} \\ 
  {} \\ 
  P' 
\end{array}} \right] = Q\left( { - V \cdot \spacevec n} \right) + \left(V \cdot \spacevec n\right)Q = 0,
\end{equation}
which is neutrally stable, having no additional stabilizing dissipation. We note again, that the mean state for the linearization is chosen such that the normal mean velocity components are zero, resulting in the zeros on the diagonal of $\mmatrixAn$.

Substituting the exact upwind flux where $\varepsilon=1$ into \eqref{eq:LinearStabilityCondition} and rearranging,
\begin{equation}
B_{L,1}= - \statevec U^{T}\left|  {\mmatrixAn^{-}}\right| \statevec U^{\ext}+\oneHalf\statevec U^{T} \left|  {\mmatrixAn}\right|\statevec U,
\label{eq:BL1}
\end{equation}
where $\mmatrixAn^{-}=\oneHalf\left(\mmatrixAn - \left|\mmatrixAn\right|\right)$ is negative semidefinite. 
The second term is non-negative, depends only on the interior state, and adds stabilizing dissipation. From the matrix absolute value, the dissipation term is
\begin{equation}
\statevec U^{T} \left|  {\mmatrixAn} \right|\statevec U = \frac{1}{\bar c}Q^{2}+\bar c^{3}\Ma^{2}_{n},
\end{equation}
where $\Ma_n=V_n/\bar{c}$ is the normal Mach number. Stability depends, then, on the value of the first term, which is where the boundary conditions are incorporated through the external state $\statevec U^{\ext}$ written in \eqref{eq:LinearExtState}. Then
\begin{equation}
\statevec U^{T}\left|  {\mmatrixAn^{-}}\right| \statevec U^{\ext} = \frac{1}{2\bar c}Q^{2} -\frac{\bar c^{3}}{2}\Ma^{2}.
\end{equation}
Therefore, using the upwind numerical flux, \eqref{eq:BL1} becomes
\begin{equation}
B_{L,1}=\bar c^{3}\textrm{Ma}^{2}_{n}\ge 0,
\end{equation}
as required. The amount of dissipation depends on how far the interior computed normal velocity deviates from zero.

The combination of the reflective state and local Lax-Friedrichs flux is also linearly stable. In that case the exact matrix absolute value is replaced by a diagonal matrix, $\left|\mmatrix A_{n}\right|\approx \left|\lambda\right|_{\mathrm{max}}\mmatrix I$. The jump term is added to the central (averaged) flux so
\begin{equation}
B_{L,LF} =-\frac{\left|\lambda\right|_{\mathrm{max}}}{2}\statevec U^{T}\left(\statevec U^{\ext}-\statevec U\right)=\bar c^{2}{\left|\lambda\right|_{\mathrm{max}}}\Ma_{n}^{2}\ge 0
\end{equation}

Finally, a dissipative version of the direct numerical flux \eqref{eq:NormalStarFlux} can be formed by looking at the reflective state approach. For instance, the equivalent to using the Lax-Friedrichs flux is to choose $\rho^{*}=\rho'$ and
\begin{equation}
P^{*} = P' + \frac{\bar c^{3}}{a}{\left|\lambda\right|_{\mathrm{max}}}\Ma_{n}.
\end{equation}
 Then $Q^{*} = Q + \bar c^{3}{\left|\lambda\right|_{\mathrm{max}}}\Ma_{n}$ and 
\begin{equation}
V_{n}\left( Q^{*}-Q\right) = \bar c^{2}{\left|\lambda\right|_{\mathrm{max}}}\Ma_{n}^{2}\ge 0.
\end{equation}
A similar, though more complicated, modified $P^{*}$ can be made to be equivalent to the exact upwind flux.

\subsection{Entropy Stability of Wall Boundary Condition Approximations}
As in the linear approximation, the wall boundary condition can be imposed for the nonlinear equations either by directly specifying the numerical flux or by computing it through a Riemann solver using a reflection external state that enforces the normal wall condition implicitly.
Note that in this section, the discrete variables $(\rho,\spacevec{V},P)$ describe the full nonlinear state.

For the nonlinear equations, we construct the numerical flux for a slip-wall as
\begin{equation}
\left(\bigstatevec F \cdot \spacevec n\right)^*  = \left[ {\begin{array}{*{20}c}
   0  \\
    P^*\,\spacevec n  \\
   0  \\
 \end{array} } \right]
\end{equation}
where we imposed $V_n=0$ leading to a flux with no mass or energy transfer, and we introduce a wall pressure $P^*$, whose value will be chosen to ensure consistency and stability.

%In the discrete case, the  entropy inequality condition at boundary faces [BR1paper] $\partial E \subseteq \partial \Omega$   reads as
%\begin{equation}
% \statevec w^T \left(\left(\bigstatevec f\cdot \spacevec n\right)^* 
% -\bigstatevec f\left(\left.\statevec u\right|_{\partial E}\right)\cdot\spacevec n\right)
% +\spacevec f^\ent\left(\left.\statevec u\right|_{\partial E}\right)\cdot\spacevec n \geq 0   
%\end{equation} 
After some manipulations, the discrete entropy stability condition \eqref{eq:EntropyStabilityCondition} becomes
\begin{equation}
 \begin{aligned}
   -\rho  {V_n}\left(\frac{\gamma-\varsigma}{(\gamma-1)}-\beta|| {\spacevec V}||^2\right) 
   + 2\beta V_n \left(P^*- {P}-\rho || {\spacevec V}||^2\right)
   + 2\beta V_n \left(\rho E+P\right)
   - \frac{\rho \varsigma V_n}{(\gamma-1)}=&\\
   -\rho V_n\left(\frac{\gamma}{(\gamma-1)} -\beta|| {\spacevec V}||^2\right)
   + 2\beta V_n\left(P^*+\rho E-\rho|| {\spacevec V}||^2\right)=&\\
   \frac{\rho V_n}{P}\left(-\frac{\gamma}{(\gamma-1)}P +\frac{1}{2}\rho|| {\spacevec V}||^2 
   +   P^*+\frac{P}{(\gamma-1)}-\frac{1}{2}\rho|| {\spacevec V}||^2\right)   =&\\
   \rho V_n\left(\frac{P^*}{P}-1\right)&\geq 0
   \label{eq:entropy_ineq_wall}
 \end{aligned}
\end{equation}
Therefore if we choose $P^*=P$, to be the internal pressure, the boundary flux does not contribute to the total entropy, independent of the inner normal velocity $V_n$. A value of $P^{*}$ that leads to a dissipative boundary condition can be found either through exact solution of the Riemann problem at the boundary, or through the use of an external state and an approximate Riemann solver.

\subsubsection{Exact solution of the Riemann problem}
In \cite{VanDerVegt2002_BC} a symmetric 1D Riemann problem is exactly solved following Toro~\cite{torobook}, to get the wall pressure $P^*$, accounting for the fact that $V_n$ never vanishes discretely and therefore the wall pressure should be different from the interior pressure.
The exact solution of the 1D Riemann problem reads as
\begin{equation}
 \left(\frac{P^*}{P}\right)_\text{RP}=\left\{\begin{array}{*{20}lll}
    1+\gamma\Ma_n\left(\frac{(\gamma+1)}{4}\Ma_n+\sqrt{\left(\frac{(\gamma+1)}{4}\Ma_n\right)^2+1}\right)&> 1 &\quad\text{for}\quad V_n>0 \\%\quad\text{\emph{(shock)}} \\
   \left(1+\frac{1}{2}(\gamma-1)\Ma_n\right)^{\frac{2\gamma}{(\gamma-1)}}&\leq 1 &\quad\text{for}\quad  V_n\leq 0 %\quad\text{\emph{(rarefaction)}} \\
 \end{array}\right.\label{eq:pstarVDV}
\end{equation} 
with the normal Mach number, $\Ma_n=\frac{V_n}{c}$, and the sound speed $c=\sqrt{\gamma \frac{P}{\rho}}$. 

As shown by Toro~\cite{torobook}, the solution for the rarefaction has a limiting vacuum solution for $\Ma_n \leq -2(\gamma-1)^{-1}.$ We will restrict our analysis to normal Mach numbers yielding strictly positive pressure solutions only  ($\Ma_n>\,-5$ for $\gamma=\frac{7}{5}$).

%Note that in the case of the rarefaction, if $\Ma_n <-\frac{4}{\gamma-1}$ and the exponent $\frac{2\gamma}{\gamma-1}$ is an even number, $\frac{P^*}{P}\leq1$ will not hold anymore. For $\gamma=7/5$ or $\gamma=5/3$, the exponent is odd and therefore $\frac{P^*}{P}\leq1$ also holds for the rarefaction. The rarefaction would not be even for $\gamma=\frac{y}{y-1}\;\;\forall\;\; y\in\mathbb{N}$.  
%If Molecular theory is stricly appield to an ideal gas, with $M$ being the number of degrees of freedom of the molecule, $\gamma=\frac{M+2}{M}$ and the exponent becomes $\frac{2\gamma}{\gamma-1}=M+2$, with $3\leqM\leq\infty$   

It is easy to see that using $P^*$ from \cref{eq:pstarVDV}, the entropy inequality \cref{eq:entropy_ineq_wall} is still satisfied for $|V_n| \neq 0$, and the added entropy scales with the discrete value of $V_n$ at the boundary. Hence, for $h\rightarrow0$, the discrete boundary condition converges to its physical counterpart, since $V_n\rightarrow0$. The choice of $P^*$ from \cref{eq:pstarVDV} appears to stabilize under-resolved simulations, which can be now explained by the fact that the boundary flux always adds entropy for $|V_n| \neq 0$. 

\subsubsection{Using approximate Riemann solvers for the boundary flux}
A well known strategy  in finite volume methods is to mirror only the velocity of the internal state and solve an approximate Riemann problem to get the boundary flux, mostly just because of a simpler implementation, since an approximate Riemann solver is already available and used for the fluxes between the elements. For DG methods, see also, for example, \cite{CHEN2017427} and  \cite{PARSANI201588} where reflection conditions are proved to be entropy stable.

The mirror state is set so that the mass and energy flux are zero. Let the inner state be labeled $L$ and the outer $R$. then the inner and outer states that satisfy the mirror condition are
\begin{equation}
\label{eq:uLuR}
 \statevec U^L = \left[\rho\;\;\rho\spacevec V_n\;\; E\right]^T\,,\quad \statevec U^R=\left[\rho\;\;\rho(\spacevec V - 2V_n\spacevec n )\;\; E\right]^T
\end{equation} 
We show below under what conditions on the normal velocity $V_{n}$ that the reflection condition is entropy stable for the Lax-Friedrichs, HLL and HLLC, Roe and  {EC-ES} fluxes.

\paragraph{\textbf{Lax-Friedrichs Flux}}
We start with the simplest approximate Riemann solver, the Lax-Friedrichs or Rusanov flux, which reads as
\begin{equation}
 \left(\bigstatevec F \cdot \spacevec n\right)_\text{LF}^*=\frac{1}{2}\spacevec n\cdot\left(\bigstatevec F (\statevec U^L)+\bigstatevec F (\statevec U^R)\right)-\frac{\left|\lambda\right|_\mathrm{max}}{2} (\statevec U^R-\statevec U^L).
\end{equation} 
Inserting the states from \cref{eq:uLuR}, we get
\begin{equation}
\left(\bigstatevec F \cdot \spacevec n\right)_\text{LF}^* = \left[ {\begin{array}{*{20}c}
   {0}  \\
   {(\rho V_n^2 + P)\,\spacevec n}  \\
   {0}  \\
 \end{array} } \right]-\frac{\lambda_\mathrm{max}}{2}  \left[ {\begin{array}{*{20}c}
   {0}  \\
   {-2\rho V_n\spacevec n}  \\
   {0}  \\
 \end{array} } \right] = \left[ {\begin{array}{*{20}c}
   {0}  \\
   {(\rho V_n^2+\rho V_n\lambda_\mathrm{max} + P)\,\spacevec n}  \\
   {0}  \\
 \end{array} } \right].
\end{equation}
The maximum wave speed is normally approximated from the largest leftgoing and rightgoing wave speed, 
\begin{equation}
 \lambda_\mathrm{max}=\max(|V_n^L|+c^L,|V_n^R|+c^R)=|V_n|+c \,,\quad \text{since}\quad c^L=c^R=c\,,\quad V_n=V_n^L=-V_n^R
\end{equation} 
and thus gives a definition of $P^*$
\begin{equation}
 \left(\frac{P^*}{P}\right)_\text{LF} = 1+\gamma\Ma_n\left(\Ma_n+|\Ma_n|+1\right)=\left\{\begin{array}{*{20}lll}
    1+\gamma\Ma_n(2\Ma_n+1)&> 1 &\quad\text{for}\quad V_n>0  \\
    1+\gamma\Ma_n &\leq 1 &\quad\text{for}\quad  V_n\leq 0  \\
 \end{array}\right.,
 \label{eq:pstarLF}
\end{equation}  
which shows that the Lax-Friedrichs flux satisfies the entropy inequality \cref{eq:entropy_ineq_wall}.

\paragraph{\textbf{HLL and HLLC Flux}}
The HLL flux \cite{torobook} is written as
 \begin{equation}
 \left(\bigstatevec F \cdot \spacevec n\right)_\text{HLL}^*=\frac{1}{S^R-S^L}\left(\spacevec n\cdot\left(S^R\bigstatevec F (\statevec U^L)-S^L\bigstatevec F (\statevec U^R)\right)+S^L S^R\left(\statevec U^R-\statevec U^L\right)\right).
\end{equation} 
The leftgoing and rightgoing wave speeds are $S^L=V_n^L -c^L= -V_n^R-c^R= -S^R$ and the HLL flux reduces to
 \begin{equation}
 \left(\bigstatevec F \cdot \spacevec n\right)_\text{HLL}^*=\frac{1}{2} \spacevec n\cdot\left(\bigstatevec F (\statevec U^L)+\bigstatevec F (\statevec U^R)\right) -\frac{S^R}{2}\left(\statevec U^R-\statevec U^L\right).
\end{equation} 
If we would choose $S^R$ to be the maximum wave speed, the HLL flux would reduce to the Lax-Friedrichs flux. However, with $S^R=V_n^R+c^R =-V_n+c$, an even simpler relation for $P^*$ is found, which also satisfies the entropy inequality
\begin{equation}
 \left(\frac{P^*}{P}\right)_\text{HLL} = 1+\gamma\Ma_n\left\{\begin{array}{*{20}ll}
    > 1 &\quad\text{for}\quad V_n>0  \\
    \leq 1 &\quad\text{for}\quad  V_n\leq 0  \\
 \end{array}\right.\label{eq:pstarHLL}
\end{equation}  
For the HLLC flux \cite{torobook}, one can show that since the Riemann problem is symmetric, the approximate wave speed of the contact discontinuity is $\lambda^*=0$ and, choosing $S^R=-V_n+c$,  HLLC reduces to the HLL flux.
\begin{equation}
 \left(\frac{P^*}{P}\right)_\text{HLLC} =1+\gamma\Ma_n=\left(\frac{P^*}{P}\right)_\text{HLL} \label{eq:pstarHLLC}
\end{equation} 

\paragraph{\textbf{Roe Flux}}
For the original Roe method without entropy fix \cite{torobook}, the mean values are
\begin{equation}
 \tilde{V}_n = \frac{\sqrt{\rho^L}V_n^L +\sqrt{\rho^R}V_n^R}{\sqrt{\rho^L} +\sqrt{\rho^R}}=0\,,\quad \tilde{V}_{t_1}=V_{t_1}\,,\quad \tilde{V}_{t_2}=V_{t_2}\,,\qquad
 \tilde{c} = c\sqrt{1+\frac{(\gamma-1)}{2}\Ma_n^2} .
\end{equation} 
After some manipulations,  
\begin{equation}
\left(\bigstatevec F \cdot \spacevec n\right)_\text{Roe}^* = \left(\bigstatevec F \cdot \spacevec n\right) + \tilde{\lambda}_1\tilde{\alpha}_1\tilde{\statevec K}^{1}= \left(\bigstatevec F\cdot \spacevec n\right)+(-\tilde{c})\frac{\rho V_n}{\tilde{c}}\left[ {\begin{array}{*{20}c}
   {1}  \\
   {-\tilde{c}}  \\ V_{t_1}\\V_{t_2} \\
   { \frac{1}{\rho}(\rho E + P)}  \\
 \end{array} } \right]=\left[ {\begin{array}{*{20}c}
   {0}  \\
   {(\rho V_n^2 +\rho V_n \tilde{c} + P)\,\spacevec n}  \\
   {0}  \\
 \end{array} } \right]\,.
\end{equation}
with $\tilde{\lambda_1}=\tilde{V}_n-\tilde{c}=-\tilde{c}$, $\alpha_1=\rho V_n /\tilde{c}$ and  $\tilde{\statevec K}^{1}$ from \cite{torobook}. 
This leads again to a definition of $P^*$ 
\begin{equation}
 \left(\frac{P^*}{P}\right)_\text{Roe} =1+\gamma \Ma_n\left(\Ma_n+\sqrt{1+\frac{(\gamma-1)}{2}\Ma_n^2}\right),
 \label{eq:pstarRoe}
\end{equation} 
which fulfills the entropy inequality as long as 
\begin{equation}
 \Ma_n \geq -\sqrt{\frac{2}{3-\gamma}} \,,\quad \text{for}\,\gamma=\frac{7}{5}\quad \Ma_n > -1.12 \,.
\end{equation} 
Thus, the Roe flux is entropy stable for shocks, but not for supersonic rarefactions. 
%This result is expected, as it is well known that the original Roe method violates the entropy condition  \cite{torobook}.

\paragraph{\textbf{EC-ES Fluxes}}
We can also apply an entropy conservative (EC) flux that is used for interior element interfaces and add an entropy stable dissipation term (ES) to compute the boundary flux via the mirrored states \cref{eq:uLuR}. This is exactly the strategy proposed in Parsani et al. \cite{PARSANI201588} to get the boundary flux. Such an EC-ES flux is presented in Winters et al. \cite{winters2017_ESmatrix} %, where an entropy conserving (EC) flux is combined with an entropy stable dissipation term. The numerical flux is
 \begin{equation}
 \left(\bigstatevec F \cdot \spacevec n\right)_\text{ES}^*=\bigstatevec F_\text{EC}\left(\statevec U_L,\statevec U_R\right) \cdot \spacevec n - \frac{1}{2} \mmatrix{D}\,\mmatrix H \left(\statevec W^R-\statevec W^L\right) \label{eq:ESflux}
\end{equation} 
where $\mmatrix{D}$ is a dissipation matrix and the matrix $\mmatrix{H} \jump{\statevec w} \simeq \jump{\statevec{u}}$ is carefully derived from the left and right states. Details are given in \cite{winters2017_ESmatrix}, where two approaches for the dissipation are distinguished. One is a Lax-Friedrichs-type dissipation, scaling with the maximum eigenvalue $\lambda_\mathrm{max}=|V_n|+c$ (referred to as 'EC-LF'). The other is a Roe-type dissipation computed via the eigenstructure of the matrix $(\mmatrix{D}\,\mmatrix{H})$ (referred to as 'EC-Roe').

If we carefully insert the two mirrored boundary states into \cref{eq:ESflux}, we again get an equation for the modified pressure
\begin{equation}
 \left(\frac{P^*}{P}\right)_\text{EC-LF} =1+\gamma \Ma_n\left(|\Ma_n|+1\right)
 \label{eq:pstarEC-LF}
\end{equation} 
for the Lax-Friedrichs-type dissipation and 
\begin{equation}
 \left(\frac{P^*}{P}\right)_\text{EC-Roe} =1+\gamma \Ma_n 
 \label{eq:pstarES}
\end{equation} 
for the Roe-type dissipation. Both approaches lead to an entropy stable boundary flux when using a mirrored state. Note that the modified pressure of the EC-Roe flux \eqref{eq:pstarES} exactly matches the one of the HLL flux \cref{eq:pstarHLL}.

% The entropy fix by Harten-Hyman  introduces a modified wavespeed $\bar{\lambda}_1$, so that
% \begin{equation}
% \left(\bigstatevec f \cdot \spacevec n\right)_{\text{Roe},\epsilon}^* = \left(\bigstatevec f\cdot \spacevec n\right) + \bar{\lambda}_1\tilde{\alpha}_1\tilde{\statevec K}^{(1)}= \left(\bigstatevec f\cdot \spacevec n\right)+\bar{\lambda}_1\frac{\rho V_n}{\tilde{c}}\left[ {\begin{array}{*{20}c}
%    {1}  \\
%    {-\tilde{c}}  \\ v_{t_1}\\v_{t_2} \\
%    { \frac{1}{\rho}(\rho E + p)}  \\
%  \end{array} } \right]= 
%  \left[ {\begin{array}{*{20}c}
%    {\rho V_n\left(1+\frac{\bar{\lambda}_1}{\tilde{c}}\right)}  \\
%    {(\rho V_n \bar{\lambda}_1 + p)\,\spacevec n+\rho\spacevec v V_n \left(1+\frac{\bar{\lambda}_1}{\tilde{c}}\right)}   \\
%    {(\rho E + p)V_n\left(1+\frac{\bar{\lambda}_1}{\tilde{c}}\right)}  \\
%  \end{array} } \right]
% \end{equation}
% which does not preserve mass or energy and should not be used as boundary condition flux...\\
% ...other E-fix?

% \begin{figure}[!htb]
% \centering\includegraphics[width=0.98\textwidth]{pstar_p.pdf}\\
% \includegraphics[width=0.98\textwidth]{pstar_p_loglog.pdf}\\
% \caption{\label{fig:pstar_p} Ratio $\frac{P^*}{P}$ and entropy contribution $\Delta s$ over the normal Mach number $\Ma_n$, with a log-log scale in the lower plot, for $\gamma=1.4$. RP refers to \cref{eq:pstarVDV}, LF to \cref{eq:pstarLF}.}
% \end{figure}

\section{Discussion}
In the previous section we have shown conditions under which a specified wall flux is stable. In the linear analysis, the central numerical flux adds no dissipation and is neutrally stable. In the nonlinear analysis, entropy is not generated if the numerical wall pressure is equal to the internal pressure, $P^{*}=P'$. For upwinded approximations, the amount of energy or entropy dissipation depends on the normal Mach number. Since the boundary condition is only imposed weakly through the numerical flux, the normal Mach number will not be exactly zero except in the convergence limit. In fact, flow computations (especially steady state ones) are usually initiated with an impulsive start, where the initial state is a uniform flow, and the normal Mach number is not zero. This has proved over time to be very robust in practice. The analysis above gives an explanation why.

In the linear analysis the dissipation due to imposing the boundary condition is proportional to the square of the normal Mach number. With an impulsive start initialization, this dissipation will be large. As the flow develops and the boundary condition is better enforced, the dissipation reduces, going away only as the approximate solution converges.

\begin{figure}[!htbp]
\centering
\includegraphics[width=0.85\textwidth]{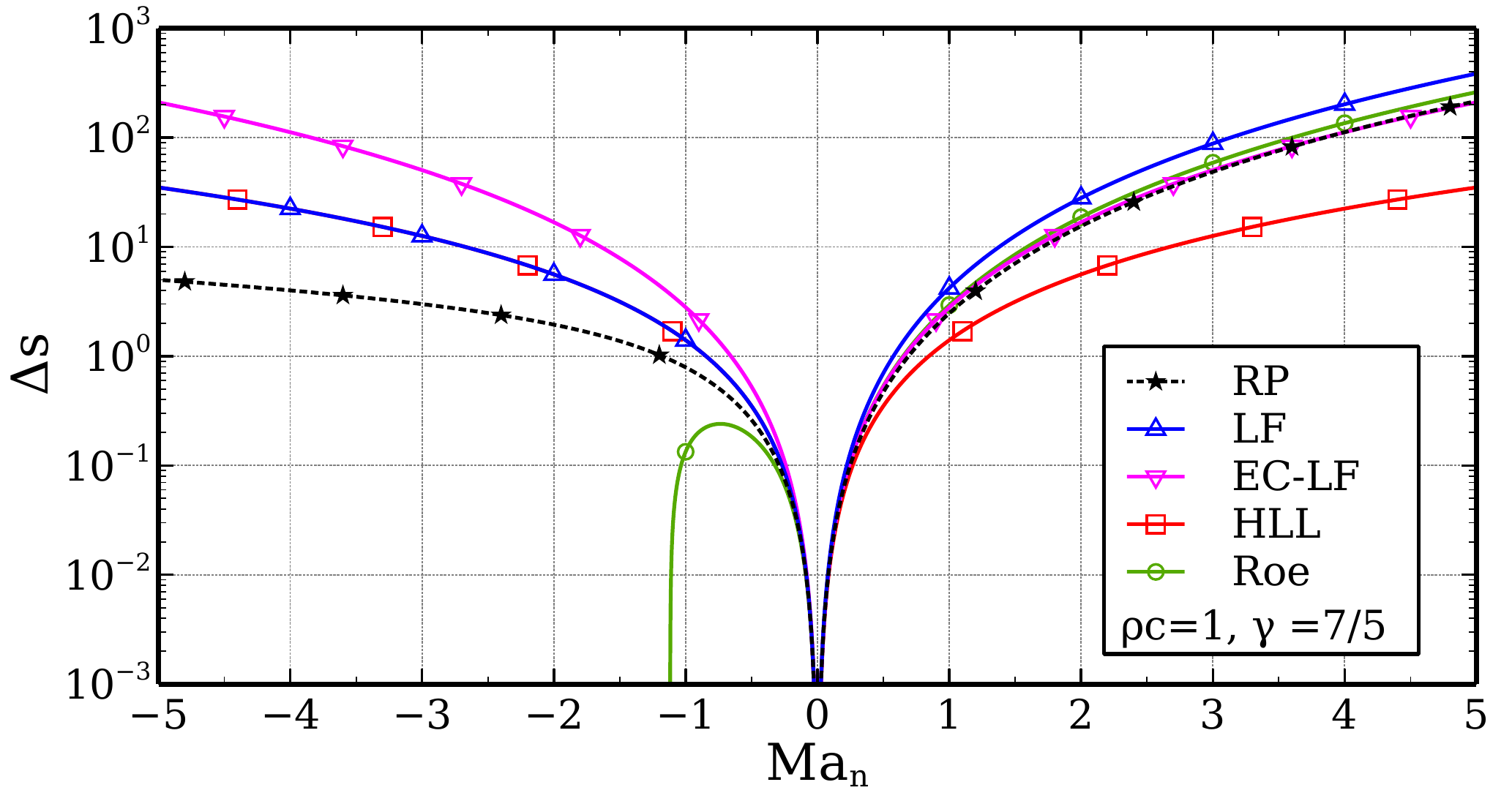}\\
\includegraphics[width=0.85\textwidth]{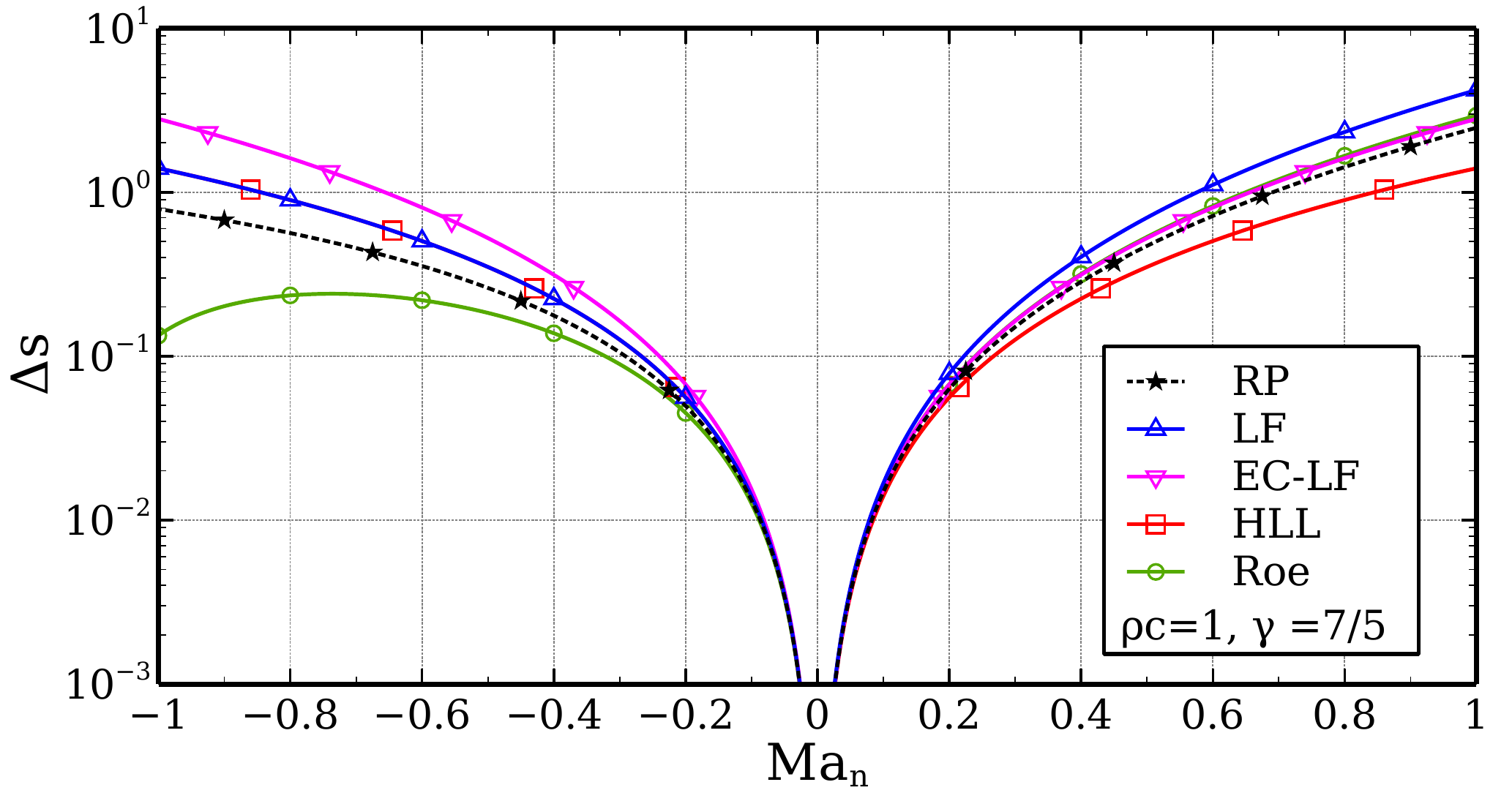}\\
\caption{\label{fig:delta_s} Entropy contribution $\Delta s$ \cref{eq:delta_s} produced by the wall boundary flux. RP refers to the exact Riemann problem \cref{eq:pstarVDV}, LF to \cref{eq:pstarLF}, EC-LF to \cref{eq:pstarEC-LF}, HLL to \cref{eq:pstarHLL} and Roe to \cref{eq:pstarRoe}. Plotted over the normal Mach number ranges $|\Ma_n|\leq 5$ on the top and restricted to $|\Ma_n|\leq 1$ on the bottom.}
\end{figure}

A similar effect is observed for the use of the different approximate Riemann solvers in the nonlinear analysis. In \cref{fig:delta_s}, we compare the entropy contribution
\begin{equation}
 \Delta s = (\rho c) \Ma_n\left(\frac{P^*}{P}-1\right) \label{eq:delta_s}
\end{equation} 
for the different wall boundary fluxes, over a range of normal Mach numbers for $(\rho c) =1$ and $\gamma=7/5$. When the boundary condition is exactly fulfilled ($\Ma_n=0$), the entropy contribution is zero. For low normal Mach numbers, all fluxes have the same behavior. Compared to the exact Riemann problem (RP), the Lax-Friedrichs flux and the EC-LF flux always produce more entropy whereas the the HLL flux produces less entropy for impinging velocities $\Ma_n>0$. The results of HLLC and EC-Roe fluxes are not plotted, as they coincide with the HLL flux. As shown in the analysis, the Roe flux produces a negative entropy change for supersonic rarefactions, implying that it is not suitable for all flow configurations.

\acknowledgement{ This work was supported by a grant from the Simons Foundation (\#426393, David Kopriva). G.G. has been supported by the European Research Council (ERC) under the European Union's Eights Framework Program Horizon 2020 with the research project \textit{Extreme}, ERC grant agreement no. 714487. Florian Hindenlang thanks Eric Sonnendr\"ucker and the Max-Planck Institute for Plasma Physics in Garching for their constant support. 
%We would also like all participants of the ICOSAHOM 2018 for the valuable discussions on the topic of entropy stable schemes, which motivated this work.
}

\bibliographystyle{plain}

\bibliography{dakBib}

\end{document}